\documentclass[12pt]{amsart}
\usepackage{amscd}
\usepackage{hyperref}
\usepackage[all]{xypic}
\def\depth{\operatorname{depth}}
\newcommand{\m}{\frak m}

\newtheorem{lemma}{Lemma}[section]
\newtheorem{corollary}[lemma]{Corollary}
\newtheorem{theorem}[lemma]{Theorem}
\newtheorem{proposition}[lemma]{Proposition}

\newtheorem{remark}[lemma]{Remark}

\newtheorem{example}[lemma]{Example}

\advance\headheight1.15pt
\begin{document}
\title{On the lengths of quotients of ideals and depths of fiber cones}

\author{A. V. Jayanthan}
\address{Department of Mathematics, Indian Institute of Technology
Madras, Chennai, INDIA - 600036.}
\thanks{$^\dag$ Supported by the Council of Scientific and Industrial
Research, India}
\thanks{AMS Classification 2000: 13H10, 13H15 (Primary), 13C15, 13A02
(Secondary)}
\email{jayanav@iitm.ac.in}
\author{Ramakrishna Nanduri$^\dag$}
\address{Department of Mathematics, Indian Institute of Technology
Madras, Chennai, INDIA - 600036.}
\email{ramakrishna@iitm.ac.in}

\maketitle

\begin{abstract}
Let $(R,\m)$ be a Cohen-Macaulay local ring, $I$ an $\m$-primary
ideal of $R$ and $J$ its minimal reduction. We study the depths of
$F(I)$ under certain depth assumptions on $G(I)$ and length condition
on quotients of powers of $I$ and $J$, namely $\sum_{n\geq0}\lambda(\m
I^{n+1}/\m JI^n)$ and $\sum_{n\geq0}\lambda(\m I^{n+1} \cap J/\m
JI^n)$.
\end{abstract}

\section{Introduction}
Let $(R,\m)$ be a Cohen-Macaulay local ring of dimension $d > 0$ with
infinite residue field and $I$ an $\m$-primary ideal of $R$. The fiber
cone of $I$, denoted by $F(I) := \oplus_{n\geq 0}I^n/\m I^n$, the
associated graded ring of $I$, denoted by $G(I) := \oplus_{n\geq
0}I^n/I^{n+1}$ and the Rees algebra $\mathcal{R}(I) := \oplus_{n\geq0}
I^nt^n \subset R[t]$ are together known as blowup algebras related to
$I$. Many authors have studied the relationship between properties of
the ideal and its blowup algebras.  Northcott and Rees introduced the
concept of a reduction of an ideal to study various properties of the
ideal and its blowup algebras. An ideal $J \subseteq I$ is said to be
a reduction of $I$ with respect to an $R$-module $M$ if $I^{n+1}M =
JI^nM$ for some $n \geq 0$. The integer $r^M_J(I) = \min\{n \; | \;
I^{n+1}M = JI^nM\}$ is called the $M$-reduction number of $I$ with
respect to $J$. If $M = R$, then $J$ is said to be a reduction of $I$
and the integer $r_J(I) = r^R_J(I)$ is known as the reduction number
of $I$ with respect to $J$. A reduction is said to be a minimal
reduction if it is minimal with respect to inclusion. It is known that
if the residue field of $R$ is infinite, then all minimal reductions
are minimally generated by $\ell(I)$ number of elements, where
$\ell(I) := \dim F(I)$ is the analytic spread of $I$.

The relation between the lengths of quotients of ideals and depths of
blowup algebras has been a subject of several papers. Let $\lambda(-)$
denote the length function. It has been shown by many authors that the
two integers 
$$
\Delta(I,J) = \sum_{p\geq1}\lambda\left(\frac{I^{p+1}\cap
J}{I^pJ}\right), \hspace*{1cm} \Lambda(I,J) = \sum_{p\geq0}
\lambda\left(\frac{I^{p+1}}{I^pJ}\right)
$$
controls the depth of the associated graded ring.  Valabrega and Valla
proved that $I^n \cap J = JI^{n-1}$ for all $n \geq 1$ if and only if
$G(I)$ is Cohen-Macaulay, where $J$ is a minimal reduction of $I$,
\cite{vv}. The Valabrega-Valla condition can be rephrased as
$\sum_{n\geq 1}\lambda(I^n \cap J/JI^{n-1}) = 0$. Guerrieri studied
ideals satisfying $\sum_{n\geq 1}\lambda(I^n \cap J/JI^{n-1}) = 1$ and
showed that in this case $\depth G(I) = d-1$, \cite{g}. She also
proved that if $\lambda(I^2 \cap J/JI) = 2$ and $I^k \cap J =
JI^{k-1}$ for all $k \geq 3$, then $\depth G(I) \geq d-2$, \cite{g2}.

Huckaba and Marley proved that $e_1(I) \leq \Lambda(I,J)$ and if the
equality holds, then $\depth G(I) \geq d-1$, \cite{hm}. Wang showed
that if $\Lambda(I,J) - e_1(I) = 1$, then $\depth G(I) \geq d-2$. As a
consequence he deduced that if $\Delta(I,J) = 2$, then $\depth G(I)
\geq d-2$, \cite{w}. Guerrieri and Rossi proved that if $\lambda(I^2
\cap J/JI) = 3$ and $J \cap I^n = JI^{n-1}$ for all $n \geq 3$, then
$\depth G(I) \geq d-3$, assuming that $R/I$ is Gorenstein and $d \geq
4$, \cite{gr}.

In the case of fiber cone, the similar relations have not been
investigated well. Cortadellas and Zarzuela proved that if $G(I)$ is
Cohen-Macaulay, then $F(I)$ is Cohen-Macaulay if and only if 
$\sum_{n\geq2}\lambda(\m I^n \cap J/\m I^{n-1}J) = 0$, \cite{cz}. It
is known that $r(I) \leq 1$ implies that the fiber cone is
Cohen-Macaulay, \cite{s}. The relation between $\m$-reduction number
and the depth of fiber cone is not that strong. The fiber cone need
not even have depth $d-1$ when $r^\m_J(I) = 1$, cf. Example \ref{ex1}.

In this article we study the depths of fiber cones of ideals
satisfying the properties $\sum_{n\geq1}\lambda(\m I^{n+1} \cap J/\m
JI^n) \leq 1$ and $\sum_{n\geq0}\lambda(\m I^{n+1}/\m JI^n) \leq 2$
under some depth assumptions on depth of $G(I)$. The paper is
organized in the following manner. In Section 2, we present some
preliminary lemmas needed for the proof of main theorems. In Section
3, we prove:

\vskip 2mm
\noindent
{\bf Theorem \ref{redn1}}. ~
Let $(R,\m)$ be a Cohen-Macaulay local ring of dimension $d\geq 1$
with infinite residue field, $I$ an $\m$-primary ideal
and $J\subseteq I$ a minimal reduction of $I$ such that $m I^n \cap J
= \m I^{n-1}J$ for all $n \geq 1$.  If $\depth(G(I))\geq
d-t$, then $\depth(F(I))\geq d-t+1$ for $1 \leq t \leq d$.

\vskip 2mm
\noindent
{\bf Theorem \ref{main}}. ~
Let $(R,\m)$ be a Cohen-Macaulay local ring of dimension $d\geq 2$
with infinite residue field, $I$ an $\m$-primary ideal
and $J\subseteq I$ a minimal reduction of $I$ such that $\sum_{k
\geq 2} \lambda((\m I^k \cap J)/\m I^{k-1}J)=1$.  If $\depth(G(I))\geq
d-t$, then $\depth(F(I))\geq d-t$ for $1\leq t \leq d-1$.

\vskip 2mm
\noindent
In Section 4, we prove:
\vskip 2mm
\noindent
{\bf Theorem \ref{main2}, \ref{main1}}. ~
Let $(R,\m)$ be a Cohen-Macaulay local ring of dimension $d \geq 2$
with infinite residue field. Let $I$ be any $\m$-primary ideal of $R$
and $J \subseteq I$ a minimal reduction of $I$. Suppose $\sum_{n \geq
0}\lambda(\m I^{n+1}/\m JI^n)=1 \mbox{ or }2$. If $\depth(G(I)) \geq
d-t$, then $\depth(F(I)) \geq d-t+1$, for $2 \leq t \leq d$.
\vskip 2mm
\noindent
We also study the Cohen-Macaulay property of fiber cones in these
cases.
In the last section we present some examples to support our results.
The computations have been performed in the Computational Commutative 
Algebra software, CoCoA \cite{co}.
\section{Preliminaries}

In this section we prove some technical lemmas which are required in
the proof of main theorems. Throughout this paper $(R,\m)$ denotes a
Cohen-Macaulay local ring with infinite residue field, $I$ an
$\m-$primary ideal and $J$ its minimal reduction. First, we recall some
results from the literature that we need.

\begin{theorem}\cite[Theorem 1.1.7]{bh} \label{rees}
Let $R$ be a ring, $M$ an $R-$module,$x_1,\ldots,x_n$ an $M-$regular
sequence, and $I=(x_1,\ldots,x_n)$.  Let $X_1,\ldots,X_n$  be
indeterminates over $R$. If $F \in M[X_1,\ldots,X_n]$ is homogeneous
of total degree $d$ and $F(x_1,\ldots,x_n) \in I^{d+1}M$, then the
coefficients of $F$ are in $IM$.
\end{theorem}

\begin{lemma}\cite[Lemma 5.2]{avj-jv}  \label{lem1}
Let $(R,\m)$ be a Cohen-Macaulay local ring of dimension $d > 0$ with
infinite residue field. Let $I$ be an $\m-$primary ideal and $J$ a
minimal reduction of $I$. Let $\{x_1,\ldots,x_d\}$ be a minimal
generating set for $J$ such that for some index $i$, $1 \leq i \leq d$
$\m I^n \cap (x_1,\dots, \hat{x_i} ,\ldots,x_d) \subseteq \m I^{n-1}J$
for all $n$, $1 \leq n \leq k$ for some integer $k$.  Then for all 
$1 \leq n \leq k$,
$$\m I^n
\cap (x_1,\dots, \hat{x_i} ,\ldots,x_d) = \m I^{n-1}(x_1,\dots,
\hat{x_i} ,\ldots,x_d).$$ 
\end{lemma}
The following lemma, which is very useful in detecting positive depth
property of fiber cones, is known as the ``Sally machine for fiber
cones''.
\begin{lemma}[Lemma 2.7,\cite{avj-jv}] \label{lem6}
Let $(R,\m)$ be a Cohen-Macaulay local ring, $I$ an $\m$-primary ideal
in $R$ and $x \in I$ such that $x^{*}$ is superficial in $G(I)$ and
$x^{o}$ is superficial in $F(I)$. If $\depth(F(I/(x)))\geq 1$, then
$x^{o}$ is regular in $F(I)$. 
\end{lemma}

For the rest of the section, let $(R,\m)$ be a Cohen-Macaulay local
ring of dimension $d > 0$ with infinite residue field, $I$ an
$\m$-primary ideal and $J \subseteq I$ a minimal reduction of $I$ such
that for some $k \geq 2, \; \m I^n \cap J = \m JI^{n-1}$ for $1 \leq n
< k$.
\begin{lemma} \label{lem2}
Let $\{x_1,\ldots,x_d \}$ be a minimal generating set for $J$ such
that $\m I^n \cap (x_{i_1},\ldots,x_{i_r})\subseteq \m I^{n-1}J$ for
all $n\leq k$ and for some $1 \leq i_{1} < \cdots < i_r \leq d$,
where $1 \leq r <d$. Then $\m I^n \cap (x_{i_1},\ldots,x_{i_r})= \m
I^{n-1}(x_{i_1},\ldots,x_{i_r})$ for all $n\leq k$.
\end{lemma}
\begin{proof}
Without loss of generality we may assume $i_1=1,\ldots ,i_r=r$. The
Lemma 5.2 in \cite{avj-jv} proves the case $r=d-1$. Assume $1 \leq r
\leq d-2$.  We proceed by induction on $k$. Let $k=2$. Since
$\bar{x}_1, \ldots, \bar{x}_r$ is a part of an $R/\m$-basis for $I/\m
I$, $\m I \cap (x_1, \ldots, x_r) = \m (x_1, \ldots, x_r)$.  Let
$y=a_1x_1+\cdots +a_rx_r \in \m I^2 \cap (x_1,\ldots,x_r)$. By
hypothesis $y=b_1x_1+\cdots +b_dx_d$, where $b_j\in \m I$ for all $i=
1,\ldots,d$.  Therefore $(a_1-b_1)x_1+\cdots + (a_r-b_r)x_r =
b_{r+1}x_{r+1}+ \cdots + b_dx_d$. Since $\{x_1,\ldots,x_d\}$ is a
regular sequence, $b_{r+1},\ldots,b_d \in (x_1,\ldots,x_r)\subseteq
J$. $b_{r+1},\ldots,b_d \in \m I \cap J=\m J$. Therefore
$(a_1-b_1)x_1+\cdots +(a_r-b_r)x_r \in \m J^2$. By the Theorem
\ref{rees}, $a_1-b_1,\ldots,a_r-b_r \in \m J$. This implies
$a_1,\cdots,a_r \in \m I$. Hence $y \in \m I(x_1,\ldots,x_r)$. 

\vskip 2mm
\noindent
Assume
$k \geq 3$. By induction hypothesis $\m I^n \cap
(x_{1},\ldots,x_{r})= \m I^{n-1}(x_{1},\ldots,x_{r})$ for all $n\leq
k-1$. Since $\m I^k \cap (x_1,\ldots,x_r)=\m
I^{k-1}(x_1,\ldots,x_r)+[\m I^{k-1}(x_{r+1},\ldots,x_d) \cap
(x_1,\ldots,x_r)]$ it is enough to prove:
\vskip 2mm
\noindent
{\sc Claim:} $\m I^t(x_{r+1},\ldots,x_d)^{n-t} \cap (x_1,\ldots,x_r)
\subseteq \m I^{n-1}(x_1,\ldots,x_r)$ for all integers $t,n$ such that
$0 \leq t <k$ and $t < n$. \\
We prove the claim by induction on $t$. Suppose $t=0.$ We need to
prove that $\m (x_{r+1},\ldots,x_d)^n \cap (x_1,\ldots,x_r) \subseteq
\m J^{n-1}(x_1,\ldots,x_r)$ for all $n,r \geq 1$. We prove this
statement by induction on $n$.  Since $(x_1, \ldots, x_r)$ is a
regular sequence, the case $n = 1$ is obvious. Assume that $n \geq 2$
and that the statement is true for all $h<n$. Let $s \in \m
(x_{r+1},\ldots,x_d)^n \cap (x_1,\ldots,x_r)$.
By induction hypothesis 
$$s \in \m J^{n-2}(x_1,\ldots,x_r)= \m(x_1,\ldots,x_r)^{n-1}+\cdots +
\m(x_1,\ldots,x_r)(x_{r+1},\ldots,x_d)^{n-2} \subseteq \m J^{n-1}.$$
Thus we can write
$$\sum_{|\sigma| =n-1}k_{\sigma}x_1^{\sigma_1}\cdots x_r^{\sigma_r}+
\cdots +\sum_{|\rho|=n-2}k_\rho x_rx_{r+1}^{\rho_{r+1}}\cdots
x_d^{\rho_d}=s=
\sum_{|\omega|=n}j_{\omega}x_{r+1}^{\omega_{r+1}}\cdots
x_d^{\omega_d},$$
where $k_{\sigma},\ldots,k_{\rho} \in \m$. Therefore $s$ is a homogeneous
polynomial in $x_1,\ldots,x_d$ of degree $n-1$ with coefficients from
$\m$ such that $s \in \m J^n$ and $x_1,\ldots,x_d $ is a regular
sequence. By the Theorem \ref{rees} all the coefficients, $k_\sigma,
\ldots, k_\rho$ are in $\m J$. Hence 
$$
s \in \m J(x_1,\ldots,x_r)^{n-1}+\cdots +\m J(x_1,\ldots,x_r)(x_{r+1},
\ldots, x_d)^{n-2} \subseteq \m J^{n-1}(x_1,\ldots,x_r).
$$
This proves case $t=0$. Assume $t \geq 1$. Let $s \in \m
I^t(x_{r+1},\ldots,x_d)^{n-t} \cap (x_1,\ldots,x_r)$. Then by what we
have shown in proving the case $t=0$, 
$$\m J^{n-t-1}(x_1,\ldots,x_r)=\m (x_1,\ldots,x_r)^{n-t}+\cdots +\m
(x_1,\ldots,x_r)(x_{r+1},\ldots,x_d)^{n-t-1}.$$
It is possible to write
$$\sum_{|\sigma| =n-t}k_{\sigma}x_1^{\sigma_1}\cdots x_r^{\sigma_r}+
\cdots + \sum_{|\rho|=n-t-1}k_\rho x_rx_{r+1}^{\rho_{r+1}} \cdots
x_d^{\rho_d} = s = \sum_{|\omega|=n-t}j_{\omega}x_{r+1}^{\omega_{r+1}}
\cdots x_d^{\omega_d},$$
where $k_{\sigma},\ldots,k_p \in \m$ and all $j_{\omega} \in \m I^t$.
Again using the Theorem \ref{rees}, we get $j_{\omega} \in
(x_1,\ldots,x_r)$.  Hence $j_{\omega} \in \m I^t \cap J$ for all
$\omega$. Since $t < k$, by hypothesis we obtain all $j_{\omega}
\in \m I^{t-1}J$. Therefore 
$$
s \in \m I^{t-1}J(x_{r+1}, \ldots, x_d)^{n-t} = \m I^{t-1}(x_1,
\ldots, x_r)(x_{r+1}, \ldots, x_d)^{n-t} + \m I^{t-1}(x_{r+1}, \ldots,
x_d)^{n-t+1}.$$
In particular we have shown that
\begin{eqnarray*}
\m I^t(x_{r+1},\ldots,x_d)^{n-t} \cap (x_1,\ldots,x_r) & \subseteq &
\m I^{t-1}(x_1,\ldots,x_r)(x_{r+1},\ldots,x_d)^{n-t} \\
&  & + \m I^{t-1}(x_{r+1},\ldots,x_d)^{n-t+1} \cap (x_1,\ldots,x_r).
\end{eqnarray*}
By induction hypothesis, $\m I^{t-1}(x_{r+1},\ldots,x_d)^{n-t+1} \cap
(x_1,\ldots,x_r) \subseteq \m I^{n-1}(x_1,\ldots,x_r)$. Therefore 
$$\m I^t(x_{r+1},\ldots,x_d)^{n-t} \cap (x_1,\ldots,x_r) \subseteq 
\m I^{n-1}(x_1,\ldots,x_r).$$ This completes the proof of the claim 
and hence the lemma.
\end{proof}
\begin{lemma} \label{lem3}
Let $y_1,\ldots,y_r$ be elements in $J- \m J$ such that the sets
$\{y_i,y_j\}$  $1 \leq i < j \leq r$ are part of minimal generating
sets for $J$. Then for each $0 \leq t \leq k-1$, we have $$(\m
I^t+(y_1)) \cap \cdots \cap (\m I^t+(y_r)) \subseteq \m I^t+(y_1\cdots
y_r).$$
\end{lemma}
\begin{proof}
We prove by induction on $r$. Assume $r=2$. Let $0 \leq t \leq k-1$. 
Let $w \in (\m I^t+(y_1)) \cap (\m I^t+(y_2))$. Then $i_1+w_1y_1 = w =
i_2+w_2y_2$ for some $i_1,i_2 \in \m I^t$ and $w_1, w_2 \in R$.
Therefore $w_1y_1 - w_2y_2 \in \m I^t \cap (y_1, y_2) = \m I^{t-1}
(y_1, y_2)$, by Lemma \ref{lem2}. This implies that $w_1y_1 \in
\m I^{t-1}y_1+(y_2)$. From hypothesis $\{y_1,y_2\}$ is a regular
sequence we get $w_1 \in \m I^{t-1}+(y_2)$. Hence $w \in \m
I^t+(y_1y_2)$. 

\vskip 2mm
\noindent
Assume $r \geq 3$ and that the statement is true for each integer
strictly less than $r$.  Let $w \in (\m I^t+(y_1)) \cap \cdots \cap
(\m I^t+(y_r))$. Then $w=i_r+w_ry_r$ with $i_r \in \m I^t$. Thus
$w_ry_r \in \cap_{i=1}^{r-1}(\m I^t+(y_i))$. By induction hypothesis
$w_ry_r \in \m I^t+(y_1\cdots y_{r-1})$. Therefore $w_ry_r = i +
\alpha y_1\ldots y_{r-1}$ for some $i \in \m I^t$ and $\alpha \in R$.
This gives $w_ry_r-\alpha y_1\cdots y_{r-1} \in \m I^t \cap (y_r,y_j)$
for all $j=1,\ldots,r-1$.  By the Lemma \ref{lem2}, $w_ry_r-\alpha
y_1\cdots y_{r-1} \in \m I^{t-1}(y_r,y_j)$ for all $j=1,\ldots,r-1$.
If $w_ry_r - \alpha y_1 \cdots y_{r-1} = a y_r + by_j$ for some $a, b
\in \m I^{t-1}$, then for all $1 \leq j \leq r-1$, $(w_r - a)y_r \in
(y_j)$ so that $(w_r - a) \in (y_j)$. Thus 
$w_r \in \cap_{j=1}^{r-1}(\m I^{t-1}+(y_j))$. By induction
hypothesis $w_r \in \m I^{t-1}+(y_1\ldots y_{r-1})$.  Hence $w \in \m
I^{t}+(y_1,\ldots,y_r)$.
\end{proof}
\begin{lemma} \label{lem4}
Let $y_1,\ldots,y_r$ be elements in $J- \m J$ such that the sets
$\{y_i,y_j\}$  are part of minimal generating sets for $J$ for $1 \leq
i < j \leq r$. Suppose that $y_1,\ldots,y_r$ satisfy the following
equalities: $\m I^k \cap (y_1)=\m I^{k-1}y_1+\m I^k \cap (y_1y_i)$ for
all $i=2,\ldots,r$. Then $\m I^k \cap (y_1)=\m I^{k-1}y_1+\m I^k \cap
(y_1\cdots y_r)$.
\end{lemma}
\begin{proof}
We proceed by induction on $r$. If $r=2,$ then there is nothing to
prove. Assume $r \geq 3$. By induction hypothesis, we have the
following two equalities:
$$
\m I^k \cap (y_1)=\m I^{k-1}y_1+\m I^k \cap (y_1\cdots y_{r-1})\mbox{
and } \m I^k \cap (y_1)=\m I^{k-1}y_1+\m I^k \cap (y_1y_r).
$$
Therefore $\m I^k \cap (y_1y_r) \subseteq \m I^{k-1}y_1+\m I^k \cap
(y_1\cdots y_{r-1})$. Suppose $\beta y_1y_r \in \m I^k \cap (y_1y_r)$.
Then $\beta y_1y_r=iy_1+\alpha y_1\cdots y_{r-1}$ for some $i \in \m
I^{k-1}$ and $\alpha \in R$. As $y_1$ is regular, $\beta y_r = i + 
\alpha y_2\cdots y_{r-1}$.  This implies $\beta y_r \in
\cap_{j=2}^{r-1} [\m I^{k-1} + (y_j)]$.  By the proof of Lemma
\ref{lem3}, we get $\beta \in \cap_{j=2}^{r-1}(\m I^{k-2}+(y_j))$ so
that by Lemma \ref{lem3} we get $\beta \in \m I^{k-2}(y_2\cdots
y_{r-1})$.  Therefore $\beta y_1y_r \in (\m I^{k-2}+(y_2\cdots
y_{r-1}))(y_1y_r) \subseteq \m I^{k-1}y_1+\m I^k \cap (y_1\cdots y_r)$
which in turn gives $\m I^k \cap (y_1) \subseteq \m I^{k-1}y_1+\m I^k
\cap (y_1\cdots y_r)$. Hence $\m I^k \cap (y_1)= \m I^{k-1}y_1+\m I^k
\cap (y_1\cdots y_r)$.
\end{proof}
\begin{lemma} \label{lem5}
Let $d \geq 2$. If $\lambda((\m I^k \cap J)/\m I^{k-1}J) \neq 0$ for
some $k \geq 2$ and $\m I^n \cap J= \m I^{n-1}J$ for all $n$, $1 \leq
n < k$, then there exists an element $x \in J-\m J$ such that 
\begin{enumerate}
  \item $x^* \in G(I)$ is superficial in $G(I)$ 
  \item $x^o \in F(I)$ is superficial in $F(I)$ and 
  \item $\lambda((\m I^k \cap J)/(\m I^{k-1}J+\m I^k \cap (x))) \neq
    0$.
\end{enumerate}
\end{lemma}

\begin{proof}
Since it is possible to choose a minimal generating set $\{x_1,
\ldots, x_d\}$ for $J$ such that $x_1^*, \ldots, x_d^*$ is a
superficial sequence in $G(I)$ and $x_1^o, \ldots, x_d^o$ is a
superficial sequence in $F(I)$, we show that the result holds for a
part of a minimal generating set of $J$.

\vskip 2mm
\noindent
We prove the lemma by induction on $k$.  Let $k=2$. Let
$\{x_1\ldots,x_d\}$ be a minimal basis for $J$. Suppose that
$$
\lambda((\m I^2 \cap J)/(\m IJ+\m I^2 \cap (x_i))) = 0 = \lambda((\m
I^2 \cap J)/(\m IJ+\m I^2 \cap (x_j)))
$$
for some $i \neq j$, $1 \leq i,j \leq d$. Then $\m I^2 \cap (x_j)
\subseteq \m IJ +\m I^2 \cap (x_i)$. Let $ax_j \in \m I^2 \cap (x_j)$.
Then $ax_j=a_1x_1+\cdots +a_dx_d +bx_i$ for some $a_1,\ldots,a_d \in
\m I$. This implies $(a-a_j)x_j-(b+a_i)x_i=a_1x_1+ \cdots
+\widehat{a_ix_i}+ \cdots +\widehat{a_jx_j}+\cdots + a_dx_d$. Since
$\{x_1,\ldots ,x_d \}$ is a regular sequence, $a-a_j,b+a_i \in I$ and
hence $a,b \in I$. Again from the same argument we get $a_l \in
(x_i,x_j) \subseteq J$ for all $l \neq i,j$. Since $\m I \cap J=\m J$,
all these coefficients are in $\m J$. Thus $(a-a_j)x_j-(b+a_i)x_i$ is
a homogeneous polynomial of degree $1$ in $x_i,x_j$ which belongs to
$\m J^2$. By the Theorem \ref{rees} we get $a-a_j,b+a_i \in \m J$.
But $a_i,a_j \in \m I$. Hence $a,b \in \m I$. So $ax_j \in \m Ix_j$.
We proved $\m I^2 \cap (x_j)=\m Ix_j$. However this yields the
contradiction 
$$0 = \lambda((\m I^2 \cap J)/(\m IJ+\m I^2 \cap
(x_j))) = \lambda((\m I^2 \cap J)/\m IJ) \neq 0.$$ 
Hence for each set of minimal generators for $J$ there are at least
$d-1$ elements which satisfy our requirement. 
\vskip 2mm
\noindent
Let $k \geq 3$. Let
$\{x_1,\ldots,x_d\}$ be a minimal generating set for $J$. Suppose
$\lambda((\m I^k \cap J)/(\m I^{k-1}J+\m I^{k} \cap (w)))=0$ for all
$w \in J-\m J$. In particular we have that $\lambda((\m I^k \cap
J)/(\m I^{k-1}J+\m I^{k} \cap (x_i)))=0$ for all $i=1,\ldots,d$. This
means $\m I^k \cap J= \m I^{k-1}J+\m I^k \cap (x_i)$ for all
$i=1,\ldots,d$. Fix an integer $i$ such that $2 \leq i \leq d$. Since
$\m I^k \cap (x_1) \subseteq \m I^k \cap J=\m I^{k-1}J+\m I^k \cap
(x_i)$, it can easily be seen that
$$\m I^k \cap (x_1)= \m I^{k-1}(x_1)+[\m
I^{k-1}(x_2, \cdots, \hat{x}_i, \cdots ,x_d)+\m I^k \cap (x_i)]\cap
(x_1).$$
If we choose an element $j_2x_2+\cdots +\widehat{j_ix_i}+\cdots
+j_dx_d+r_ix_i \in [\m I^{k-1}(x_2,\ldots,\hat{x}_i,\ldots,x_d)+\m I^k
\cap (x_i)]\cap (x_1)$ for some $j_h \in \m I^{k-1}$ and
$h=2,\ldots,\hat{i},\ldots,d$, then it is clear that
$j_2x_2+\cdots+\widehat{j_ix_i}+\cdots+j_dx_d \in (x_1,x_i)$.
Since $\{x_1,\ldots,x_d\}$ is a regular sequence, $j_h \in \m
I^{k-1}\cap (x_1,\ldots,\hat{x}_h,\ldots,x_d)$, for all
$h=2,\ldots,\hat{i},\ldots,d$.
By Lemma \ref{lem2}, $\m I^{k-1}\cap (x_1, \ldots,
\hat{x}_h,\ldots,x_d)= \m I^{k-2}(x_1,\ldots,\hat{x}_h,\ldots,x_d)$,
for all $h=2,\ldots ,\hat{i},\ldots,d$.  Therefore $j_h \in \m
I^{k-2}(x_1,\ldots,\hat{x_h},\ldots,x_d)$ for all $h=2,\ldots
,\hat{i},\ldots,d$ and hence $j_2x_2+\cdots +\widehat{j_ix_i}+\cdots
+j_dx_d \in \m I^{k-2}J(x_2,\ldots,\hat{x}_i,\ldots,x_d) \cap
(x_1,x_i)$.
Using the {\sc Claim} in the Lemma \ref{lem2} with $t=k-2$ and $n=k$
we can conclude that $j_2x_2+\cdots +\widehat{j_{i}x_i}+\cdots +j_dx_d \in
\m I^{k-1}(x_1,x_i)$.
Hence 
\begin{eqnarray*}
\m I^k \cap (x_1) & \subseteq & \m I^{k-1}x_1+[\m I^{k-1}(x_1,x_i)+\m I^k
\cap (x_i)]\cap (x_1)\\  & =  & \m I^{k-1}x_1+[\m I^{k-1}(x_1)+\m I^k \cap
(x_i)]\cap (x_1) \\
& = & \m I^{k-1}x_1+[\m I^{k-1}(x_1)+\m I^k \cap (x_i)\cap
(x_1)] =\m I^{k-1}x_1+\m I^k \cap (x_1x_i).
\end{eqnarray*}
The other inclusion is obvious. Therefore $\m I^k \cap (x_1)=\m
I^{k-1}x_1+\m I^k \cap (x_1x_i)$ for all $i=2,\ldots,d$. Hence, by 
Lemma \ref{lem4}  it follows that
$\m I^k \cap (x_1)=\m I^{k-1}x_1+\m I^k \cap (x_1\cdots x_d)$.
If $k < d$, then $x_1 \cdots x_d \in \m I^{k-1} x_1 \subseteq \m I^k$
so that $\m I^k \cap (x_1)=\m I^{k-1}x_1+\m I^k
\cap (x_1\cdots x_d)=\m I^{k-1}x_1$ which yields the contradiction
$$
0 = \lambda((\m I^k \cap J)/(\m I^{k-1}J+\m I^k \cap (x_1))) = 
\lambda((\m I^k \cap J)/\m I^{k-1}J) \neq 0.
$$
Suppose $k \geq d$. Since $J/\m J$ is a vector space over an infinite
field it is possible to find elements $x_h \in \cap_{i=1}^{h-1} [J-(\m
J+(x_i))]$ for $h=d+1,\ldots,k+1$ so that $\{x_1, \ldots, x_{d-1},
x_{d+1}\}, \\ \ldots, \{x_1,\ldots ,x_{d-1},x_{k+1}\}$ are minimal
generating sets for $J$. Moreover, for any $d+1 \leq h \leq k+1$, by
the selection of $x_h$, $\{\bar{x}_j, \bar{x}_h\} \in J/\m J$ is
$R/\m$-linearly independent for any $j < h$ and hence form a part of
minimal generating set for $J$. Define $y_1=x_1,\ldots,y_{d}=x_{d},\;
y_{d+1}=x_{d+1}, \ldots,y_{k+1}=x_{k+1}$. Then $\{y_i, y_j\}$ is a
part of minimal generating set of $J$ for $1 \leq i, j \leq k+1, \; i
\neq j$. Also we have $\m I^k \cap (x_1)= \m I^{k-1}x_1+\m I^{k} \cap
(x_1x_i)$ for all $i=d,\ldots,k$. Thus $y_1,\ldots,y_{k+1}$ satisfy
the hypotheses of Lemma \ref{lem4}. Therefore $\m I^k \cap (y_1)= \m
I^{k-1}y_1+\m I^{k} \cap (y_1\cdots y_{k+1})$. This implies $\m I^k
\cap (x_1)= \m I^{k-1}x_1$.
This gives the contradiction 
$$
0=\lambda((\m I^k \cap J)/(\m I^{k-1}J+\m I^k \cap (x_1)))=\lambda((\m
I^k \cap J)/\m I^{k-1}J) \neq 0.
$$
This completes the proof of the lemma.
\end{proof}

\section{Ideals with $\sum_{n\geq1}\lambda(\m I^{n+1} \cap J/\m JI^n)
\leq 1$}
In this section we study ideals satisfying the property
$\sum_{n\geq1}\lambda(\m I^{n+1} \cap J/\m JI^n) \leq 1$. It is known
that if $G(I)$ is Cohen-Macaulay, then $F(I)$ is Cohen-Macaulay if and
only if $\m I^n \cap J = \m JI^{n-1}$ for all $n \geq 1$,
\cite[Theorem 3.2]{cz}. We relax the depth condition on $G(I)$, in the
sufficiency part of the result of Cortadellas and Zarzuela. We also
give an example to show that if $\depth G(I) = d-1$ and $F(I)$ is
Cohen-Macaulay, then the equation $\m I^n \cap J = \m I^{n-1}J$ need
not hold, cf. Example \ref{ex4}.
\begin{theorem}\label{redn1}
Let $(R,\m)$ be a Cohen-Macaulay local ring of dimension $d \geq 1$
with infinite residue field, $I$ an $\m$-primary ideal and $J$ a
minimal reduction such that $\m I^n \cap J = \m I^{n-1}J$ for all $n
\geq 1$. If $\depth G(I) \geq d-t$, then $\depth F(I) \geq d-t+1$ for
all $1 \leq t \leq d$. In particular, if $\depth G(I) \geq d-1$, then
$F(I)$ is Cohen-Macaulay.
\end{theorem}

\begin{proof}
Suppose $d = 1$. Let $J = (x)$. Then $\m I^n \cap (x) = \m I^{n-1}(x)$
for all $n \geq 1$. Hence $x^o$ is regular in $F(I)$.
\vskip 2mm
\noindent
Suppose $d \geq 2$. Assume $\depth G(I) \geq d-t$. Let $x \in J$ be a
minimal generator of $J$ such that $x^*$ is superficial in $G(I)$ and
$x^o$ is superficial in $F(I)$. Let ``-'' denotes modulo $(x)$. If $t
= d$, then by induction $\depth F(\bar{I}) \geq 1$. By Lemma
\ref{lem6}, $x^o$ is regular in $F(I)$, and hence $\depth F(I) \geq
1$. If $t \leq d-1$, then $x^*$ is regular in $G(I)$ and so
$F(\bar{I}) \cong F(I)/x^oF(I)$. By induction, $\depth F(\bar{I}) \geq
d-t \geq 1$. Again by Lemma \ref{lem6}, $x^o$ is regular in $F(I)$.
Hence $\depth F(I) = \depth F(\bar{I}) + 1 \geq d-t+1$.
\end{proof}

\begin{corollary}\label{cor-redn}
\begin{enumerate}
  \item If $\m I^n \cap J = \m I^{n-1}J$ for all $n \geq 1$, then
    $\depth F(I) > 0$.
  \item If $r^\m_J(I) = 1$ and $\depth G(I)
    \geq d-t$, then $\depth F(I) \geq d-t+1$ for all $1 \leq t \leq 
    d$.
\end{enumerate}
\end{corollary}

\begin{proof}
Putting $t = d$ in Theorem \ref{redn1}, the statement (1) follows.
For (2), note that if $\m I^n = \m I^{n-1}J$, then $\m I^n \cap J =
\m I^{n-1}J$. Now the assertion follows from Theorem \ref{redn1}.
\end{proof}

Now we study the ideals satisfying the property $\sum_{n\geq
2}\lambda(\m I^n \cap J/\m I^{n-1}J) = 1$.  We prove that in this case
the depth of the fiber cone is at least as much as that of the
associated graded ring, except when the associated graded ring is
Cohen-Macaulay. We also provide an example which shows that the lower
bound on the depth of $F(I)$ is sharp.
\begin{theorem} \label{main}
Let $(R,\m)$ be a Cohen-Macaulay local ring of dimension $d\geq 2$
with infinite residue field, $I$ be an $\m$-primary ideal in $R$
and $J\subseteq I$ a minimal reduction of $I$ such that $\sum_{n
\geq 1} \lambda((\m I^{n+1} \cap J)/\m I^nJ)=1$. If
$\depth(G(I))\geq d-t$, then $\depth(F(I))\geq d-t$, for all $1\leq t 
\leq d-1$.
\end{theorem}

\begin{proof}
By hypothesis, there exists an integer $k \geq 2$ such that
$\lambda(\m I^k \cap J/\m I^{k-1}J) = 1$ and $\m I^n \cap J = \m
I^{n-1}J$ for all $n \neq k$. By Lemma \ref{lem5}, there exists $x \in
J - \m J$ such that $x^*$ is superficial in $G(I)$, $x^o$ is
superficial in $F(I)$ and $\lambda(\m I^k \cap J/\m I^{k-1}J + \m I^k
\cap (x)) = 1$.
\vskip 2mm
\noindent
We first prove the case $t = d-1$.
We do this by induction on $d$. Let $d = 2$.
Since $\lambda(\m I^k \cap J/\m I^{k-1}J) = \lambda(\m I^k \cap J/\m
I^{k-1}J + \m I^k \cap (x)) = 1$, it follows that $\m I^k \cap (x)
\subseteq \m I^{k-1}J$. Also, for all $n \neq k$, $\m I^n \cap (x)
\subseteq \m I^n \cap J = \m I^{n-1}J$. Hence we have $\m I^n \cap (x)
\subseteq \m I^{n-1}J$ for all $n \geq 1$. Since $x^*$ is superficial
in $G(I)$ and $\depth G(I) \geq 1$, $x^*$ is regular in $G(I)$.
Therefore, by Lemma \ref{lem1}, $x^o$ is a regular element in $F(I)$.
\vskip 2mm
\noindent
Assume $d \geq 3$.  Let ``-'' denotes modulo $(x)$. By repeatedly
applying Lemma \ref{lem5}, choose $\bar{x}_2, \ldots, \bar{x}_{d-1}
\in \bar{J} - \bar{\m} \bar{J}$ such that $\lambda(\bar{\m}\bar{I}^k
\cap \bar{J}/\bar{\m}\bar{I}^{k-1}\bar{J} + \bar{\m}\bar{I}^k \cap
(\bar{x}_2, \ldots, \bar{x}_{d-1})) = 1$. This implies that
$$\bar{\m}\bar{I}^{k-1}\bar{J} + \bar{\m}\bar{I}^k \cap
(\bar{x}_2, \ldots, \bar{x}_{d-1}) = \bar{\m}\bar{I}^{k-1}\bar{J}.$$
Lifting the equation back to $R$, we get
$$\m I^{k-1}J + (x) = \m I^{k-1}J + \m I^k \cap (x, x_2, \ldots,
x_{d-1}) + (x).$$
Intersecting with $\m I^k$ we get
$$
\m I^{k-1}J + \m I^k \cap (x) = \m I^{k-1}J + \m I^k \cap (x, x_2,
\ldots, x_{d-1})
$$
By the choice of $x$, $\m I^k \cap (x) \subseteq \m I^{k-1}J$.
Therefore $\m I^k \cap (x, x_2, \ldots, x_{d-1}) \subseteq \m
I^{k-1}J$. For $n \neq k$, this inequality anyway holds because of the
hypothesis. Therefore we have, $\m I^n \cap (x, x_2 \ldots, x_{d-1})
\subseteq \m I^{n-1} J$ for all $n \geq 1$. By Lemma \ref{lem1}, $\m
I^n \cap (x, x_2, \ldots, x_{d-1}) = \m I^{n-1}(x, x_2, \ldots,
x_{d-1})$. Going modulo $(x, x_2, \ldots, x_{d-2})$, we obtain
$\bar{\m}\bar{I}^n \cap (\bar{x}_{d-1}) =
\bar{\m}\bar{I}^{n-1}(\bar{x}_{d-1})$ for all $n \geq 1$. This implies
that $\bar{x}_{d-1}^o \in F(\bar{I})$ is a regular element, where ``-''
denotes modulo $(x, x_2, \ldots, x_{d-2})$, i.e., $\depth F(\bar{I})
\geq 1$. By repeatedly applying Sally machine for fiber cones, Lemma
\ref{lem6}, we get $\depth F(I) \geq 1$.

\vskip 2mm
\noindent
Now assume $d \geq 3, \; 1 \leq t \leq d-2$ and $\depth G(I) \geq
d-t$. Let $\{x_1 = x, x_2, \ldots, x_d\}$ be a minimal generating set
for $J$ such that $x_1^*, \ldots, x_d^*$ is a superficial sequence in
$G(I)$, $x_1^o, \ldots, x_d^o$ is a superficial sequence in $F(I)$ and
$\lambda(\m I^k \cap J/\m I^{k-1}J + \m I^k \cap (x_1, \ldots, x_i)) =
1$ for all $1 \leq i \leq d-t-1$. Taking modulo $(x_1, \ldots,
x_{d-t-1})$, we get $\lambda(\bar{\m}\bar{I}^k \cap
\bar{J}/\bar{\m}\bar{I}^{k-1}\bar{J}) = 1, \; \bar{\m}\bar{I}^n \cap
\bar{J} = \bar{\m}\bar{I}^{n-1}\bar{J}$ for all $n \neq k,$ and 
$\depth G(\bar{I}) \geq 1$. Therefore, by the first part of the proof,
we get $\depth F(\bar{I}) \geq 1$ so that $\bar{x}_{d-t}^o$ is regular
in $F(\bar{I})$. Since $x_1^*, \ldots, x_{d-t}^*$ is
a regular sequence in $G(I)$, $F(\bar{I}) \cong F(I)/(x_1^o, \ldots,
x_{d-t-1}^o)$. By repeated application of Lemma \ref{lem6}, we obtain
that $x_1^o, \ldots, x_{d-t}^o$ is a regular sequence in $F(I)$. Hence
$\depth F(I) \geq d-t$.
\end{proof}

\begin{remark}
By \cite[Theorem 3.2]{cz} and Theorem \ref{main}, if $\sum_{n \geq 2}
\lambda(\m I^n \cap J/\m JI^{n-1}) = 1$ and $G(I)$ Cohen-Macaulay,
then $\depth F(I) = d-1$. This has an interesting consequence, that a
necessary condition for $F(I)$ to be Cohen-Macaulay in this case is
that $\depth G(I) \leq d-1$. In Section 5, we have provided examples
with $\sum_{n\geq2}\lambda(\m I^n \cap J/\m JI^{n-1}) = 1$ and
      \begin{enumerate}
	\item $G(I)$ Cohen-Macaulay, but $F(I)$ not Cohen-Macaulay;
	\item $F(I)$ Cohen-Macaulay, but $G(I)$ not Cohen-Macaulay.
      \end{enumerate}

\end{remark}
\begin{corollary} \label{cor1}
Suppose $\lambda(\m I^2/\m IJ)=1$ and $\m I^3=\m I^2J$.  If
$\depth(G(I)) \geq d-t$, then $\depth(F(I))\geq d-t$ for all $1\leq t
\leq d-1$.
\end{corollary}
\begin{proof}
Since $(\m I^2\cap J)/\m IJ$ is a submodule of $\m I^2/\m IJ$,
$\lambda(\m I^2 \cap J/\m IJ) \leq 1$ and $\m I^n \cap J=\m I^{n-1}J$
for all $n \geq 3$. If $\lambda(\m I^2 \cap J/\m IJ)=0$, then $\m I^n
\cap J=\m I^{n-1}J$ for all $n \geq 1$. Let $\{x_1, \ldots, x_d\}$ be
a minimal generating set for $J$ such that $x_1^*, \ldots, x_d^*$ is a
superficial sequence in $G(I)$ and $x_1^o, \ldots, x_d^o$ is a
superficial sequence in $F(I)$. Then we get $\m I^n \cap (x_1, \ldots,
x_{d-t}) \subseteq \m I^{n-1}J$ for all $n \geq 1$. If $\depth G(I)
\geq d-t$, then the assertion follows from Lemma \ref{lem2} and
Theorem 28 of \cite{cz}.
\vskip 2mm
\noindent
If $\lambda(\m I^2 \cap J/\m IJ) = 1$, the assertion follows from
Theorem \ref{main}.
\end{proof}

We conclude this section by deriving a result analogues to a result by
Vasconcelos on the depth of the associated graded rings.
\begin{corollary}
Suppose $I^3=JI^2$ and $\lambda(\m I^2 \cap J/\m JI)=1$. Then
\begin{enumerate}
 \item[(a)] If $\lambda(I^2/JI)=1$, then $\depth(F(I))\geq d-1$.
 \item[(b)] If $\lambda(I^2/JI)=2$, then $\depth(F(I))\geq d-2$.
\end{enumerate}
\end{corollary}
\begin{proof}
From the hypothesis we have $\sum_{k \geq 2} \lambda((\m I^k \cap
J)/\m I^{k-1}J)=1.$ Assume $\lambda(I^2/JI)=1$ and $I^3=JI^2$. Then by
Corollary 2.3(a) in \cite{g2} we have $\depth(G(I))\geq d-1$.
Therefore by the Theorem \ref{main} we get $\depth(F(I))\geq d-1$.
This proves (a). Now assume $\lambda(I^2/JI)=2$ and $I^3=JI^2$.  Then
by the Corollary 2.3(b) in \cite{g2} we have $\depth(G(I))\geq d-2$.  Therefore by
the Theorem \ref{main} we have $\depth(F(I))\geq d-2$. This proves
(b).
\end{proof}

\section{Ideals with $\sum_{n\geq0}\lambda(\m I^{n+1}/\m JI^n) = 1
\mbox{ or } 2$}
In this section, we study fiber cones of ideals with
$\sum_{n\geq1}\lambda(\m I^n/\m I^{n-1}J) = 1$ or $2$.  In
\cite{avj-jv1}, it was proved that if $\lambda(\m I/\m J) = 1$ and
$\depth G(I) \geq d-2$, then $\depth F(I) \geq d-1$. In the following
theorem we generalize this result.

\begin{theorem} \label{main2}
Let $(R,\m)$ be a Cohen-Macaulay local ring of dimension $d\geq 2$
with infinite residue field $R/\m$. Let $I$ be an $\m$-primary ideal
in $R$ with $J\subseteq I$ a minimal reduction of $I$ such that
$\lambda(\m I/\m J)=1$. If $\depth(G(I))\geq d-t$, then 
$\depth(F(I))\geq d-t+1$ for $2 \leq t \leq d$.
\end{theorem}
\begin{proof}
We prove the result on $d$. For $d = 2$, the result follows from
Corollary 4.5 in \cite{avj-jv1}. Assume $d \geq 3$.  We first prove
the result for $t = d$. We show that if $\lambda(\m I/\m J)=1$, then
$\depth(F(I))\geq 1$.  Let $x\in J-\m
J$ be an element such that $x^{*}$ is superficial in $G(I)$ and
$x^{o}$ is superficial in $F(I)$. Let ``-'' denotes modulo $(x)$.
Since $\m I \cap (x) = \m (x), \; \bar{\m}\bar{I}/\bar{\m}\bar{J}
\cong \m I/\m J + \m I \cap (x) = \m I/\m J$.
Then $(\bar{R},\bar{\m})$ is a $d-1$ dimensional Cohen-Macaulay local
ring and $\lambda(\bar{\m} \bar{I}/\bar{\m} \bar{J})=\lambda(\m I/\m
J)=1$. Therefore by induction hypothesis $\depth(F(\bar{I}))\geq 1$.
Hence by the Lemma \ref{lem6}, $\depth(F(I))\geq 1$. This proves  the
case $t=d$. 
\vskip 2mm
\noindent
Now assume $2\leq t \leq d-1$. Choose $x\in J-\m J$ such that $x^{*}$
is superficial in $G(I)$ and $x^{o}$ is superficial in $F(I)$. Let
``-'' denotes modulo $(x)$. Then $(\bar{R},\bar{\m})$ is a $d-1$
dimensional Cohen-Macaulay local ring with $\lambda(\bar{\m}
\bar{I}/\bar{\m} \bar{J})=1$ and $\depth G(\bar{I}) \geq d-t-1$.
Induction hypothesis yields that $\depth F(\bar{I}) \geq d-t$. Since
$d-t \geq 1$, by Lemma \ref{lem6}, $x^o$ is regular in $F(I)$. Since
$x^*$ is a regular element in $G(I), \; F(\bar{I}) \cong
F(I)/x^oF(I)$. Hence $\depth F(I) = \depth F(\bar{I}) + 1 \geq d-t+1$.
\end{proof}

\begin{theorem} \label{main1}
Let $(R,\m)$ be a Cohen-Macaulay local ring of dimension $d \geq 2$.
Let $I$ be any $\m$-primary and $J \subseteq I$ a minimal reduction of
$I$ such that $\sum_{n \geq 0}\lambda(\m I^{n+1}/\m JI^n) = 2$. If
$\depth(G(I)) \geq d-t$, then $\depth(F(I)) \geq d-t+1$, for all $2
\leq t \leq d$.
\end{theorem}

\begin{proof}
First we prove the theorem for $t=d$, i.e., we show that
$\depth(F(I))\geq 1$. We do this by induction on $d$.  Suppose $d=2$.
Since $\sum_{n \geq 0}\lambda(\m I^{n+1}/\m JI^n)=2$ there are two
possible cases, namely,
\begin{enumerate}
  \item[(i)] $\lambda(\m I/\m J)=1=\lambda(\m I^2/\m JI)$ and $\m
    I^{j+1}=\m JI^j$ for all $j \geq 2$
  \item[(ii)] $\lambda(\m I/\m J)=2$ and $\m I^{j+1}=\m JI^j$ for all
    $j \geq 1.$
\end{enumerate}
In the first case, the assertion follows by Corollary 4.5 in
\cite{avj-jv1}.  Now assume that $\lambda(\m I/\m J)=2$ and $\m
I^{j+1}=\m JI^j$ for all $j \geq 1$. Let $\{x,y\}$ be a minimal
generating set for $J$ such that $x^{*},y^{*}$ is a superficial
sequence in $G(I)$ and $x^{o},y^{o}$ is a superficial sequence in
$F(I)$. To show $x^o$ is a regular element in $F(I)$ it is enough to
prove the claim below. \\
{\sc Claim :} $(\m I^{j+1}:x)=\m I^j$ for all $j\geq 0$.\\
We prove the claim by induction on $j$. Since $x$ is a part of minimal
generating set for $J$ and hence $I$, $j = 0$ case holds.  Assume
$j\geq 1$ and the induction hypothesis that $(\m I^{j}:x)=\m I^{j-1}$.
To show that $(\m I^{j+1}:x)=\m I^j$. Consider the following exact
sequence for $j \geq 1$,
$$
0 \longrightarrow (\m I^j:x)/(\m I^j:J)
\overset{\mu_y}{\longrightarrow} (\m I^{j+1}:x)/\m I^j 
\overset{\mu_x}{\longrightarrow} \m I^{j+1}/\m JI^j \longrightarrow
\bar{\m}\bar{I}^{j+1}/\bar{\m}\bar{J}\bar{I}^j \longrightarrow 0,$$
where ``-'' denotes the modulo $(x)$. The first map $\mu_y$ is
the multiplication by $y$ and second map $\mu_x$ is the multiplication
by $x$. Since $\m I^{j+1}=\m JI^j$ for all $j \geq 1$, the last two
modules of the above exact sequence are zeros. Hence the first two
modules are isomorphic. That is
$(\m I^j:x)/(\m I^j:J)\cong (\m I^{j+1}:x)/\m I^j$. Then $\m I^{j-1}
\subseteq (\m I^j:J)\subseteq (\m I^j:x)=\m I^{j-1}$, where the last
equality follows by induction hypothesis. This implies that $(\m
I^j:x)=(\m I^j:J)$. From the above isomorphism $(\m I^{j+1}:x)=\m I^j$
as required. This proves the claim. Therefore $x^{o}$ is regular in
$F(I)$. Hence $\depth(F(I))\geq 1$. 
\vskip 2mm
\noindent
Assume $d\geq 3$. Let
$\{x_1,\ldots,x_d\}$ be a minimal generating set for $J$ such that
$x_1^{*},\dots,x_d^{*}$ is a superficial sequence in $G(I)$ and
$x_1^{o},\dots,x_d^{o}$ is a superficial sequence in $F(I)$. Let ``-''
denotes modulo $(x_1)$. Then $(\bar{R},\bar{\m})$ is a
$(d-1)$-dimensional Cohen-Macaulay local ring with 
$$0 \neq \sum_{n \geq 0}\lambda(\bar{\m} \bar{I}^{n+1}/\bar{\m}
\bar{J}\bar{I}^n) \leq 2.$$
If $\lambda(\bar{\m} \bar{I}/\bar{\m} \bar{J})=1$, then the result
follows from Theorem \ref{main2} and Lemma \ref{lem6}.
If $\lambda(\bar{\m} \bar{I}^{n+1}/\bar{\m} \bar{J}\bar{I}^n)=2$, then
again the assertion follows by induction hypothesis and Lemma
\ref{lem6}. This proves the theorem for the case $t=d$. 

\vskip 2mm
\noindent
Suppose $2 \leq t \leq d-1$ and $\depth G(I) \geq d-t$.
Let $\{x_1,\ldots,x_d\}$ be a minimal generating set for $J$ such that
$x_1^{*},\dots,x_d^{*}$ is a superficial sequence in $G(I)$ and
$x_1^{o},\ldots,x_d^{o}$ is a superficial sequence in $F(I)$. 
Let ``-'' denotes the modulo $(x_1)$. Then $(\bar{R},\bar{\m})$ is a
$(d-1)$-dimensional Cohen-Macaulay local ring with
$\depth(G(\bar{I}))\geq d-1-t$ and  
$$0 \neq \sum_{n \geq 0}\lambda(\bar{\m} \bar{I}^{n+1}/\bar{\m}
\bar{J}\bar{I}^n) \leq 2.$$
If $\sum_{n \geq 0}\lambda(\bar{\m} \bar{I}^{n+1}/\bar{\m}
\bar{J}\bar{I}^n)=1$, then by Theorem \ref{main2} we get
$\depth(F(\bar{I}))\geq d-t$. If $\sum_{n \geq 0}\lambda(\bar{\m}
\bar{I}^{n+1}/\bar{\m} \bar{J}\bar{I}^n)=2$, then by induction
hypothesis, $\depth F(\bar{I}) \geq d-t$. Since $d-t \geq 1$, by Lemma
\ref{lem6}
$x_1^o$ is regular in $F(I)$. Moreover, since $\depth G(I) \geq d-t
\geq 1$, $x_1^*$ is regular in $F(I)$ so that $F(\bar{I}) \cong
F(I)/x_1^oF(I)$. Thus $\depth(F(I))=\depth(F(\bar{I}))+1 \geq d-t+1$.
\end{proof}

\begin{corollary}
If $\displaystyle{\sum_{n \geq 0}\lambda(\m I^{n+1}/\m JI^n) \leq 2}$, then
$\depth(F(I))\geq 1$.
\end{corollary}

\begin{proof}
From the proofs of Theorem \ref{main2} and Theorem \ref{main1} we get
that if $$0 \neq \sum_{n \geq 0}\lambda(\m I^{n+1}/\m JI^n) \leq 2,$$
then $\depth(F(I))\geq 1$. Now assume $\sum_{n \geq 0}\lambda(\m
I^{n+1}/\m JI^n)=0$. This implies that $\m I=\m J$. Now the result
follows from Lemma 2.3(1) of \cite{go}.
\end{proof}
\begin{corollary}
Suppose $I^3=I^2J$. Then
\begin{enumerate}
  \item[(i)] If $\lambda(I^2/IJ)=1=\lambda(\m I^2/\m JI)$, then
    $\depth(F(I))\geq d-1$.
  \item[(ii)] If $\lambda(I^2/IJ)=2$ and $\lambda(\m I^2/\m
    IJ)=1$, then $\depth(F(I))\geq d-2$.
  \item[(iii)] If $\lambda(I^2/IJ)=2=\lambda(\m I/\m J)$ and $\m
    I^2=\m IJ$, then $\depth(F(I))\geq d-2$.
\end{enumerate}
\end{corollary}
\begin{proof}
Assume $\lambda(I^2/IJ)=1=\lambda(\m I^2/\m JI)$ and $I^3=I^2J$. Then
by the Corollary 2.3(a) in \cite{g2} we get $\depth(G(I))\geq d-1$.
Hence by the Corollary \ref{cor1} we get $\depth(F(I))\geq d-1$. This
proves (i). Now assume $\lambda(I^2/IJ)=2$ and $\lambda(\m I^2/\m
JI)=1$.  Then by the Corollary 2.3(b) in \cite{g2} we have
$\depth(G(I))\geq d-2$. Hence by the Corollary \ref{cor1} we get
$\depth(F(I))\geq d-2$.  This proves (ii). Now assume
$\lambda(I^2/IJ)=2=\lambda(\m I/\m J)$ and $\m I^2=\m IJ$. Then by the
Corollary 2.3(b) in \cite{g2} we have $\depth(G(I))\geq d-2$. Since
$\depth(G(I))\geq d-2 > d-3$, by the Theorem \ref{main1} we get
$\depth(F(I))\geq d-2$. This proves (iii).
\end{proof}

\vskip 2mm
\noindent
We conclude this section by characterizing Cohen-Macaulayness of fiber
cones ideals with $\sum_{n\geq1}\lambda(\m I^n/\m JI^{n-1}) = 2$. 
\begin{proposition}
Let $(R,\m)$ be a Cohen-Macaulay local ring of dimension $d \geq 2$,
$I$ an $\m$-primary ideal such that $\depth G(I) \geq d-1$ and $J$ a
minimal reduction of $I$.
\begin{enumerate}
  \item If $\lambda(\m I/\m J) = 1 = \lambda(\m I^2/\m IJ)$ and $\m
    I^{n+1} = \m I^nJ$ for all $n \geq 2$, then
    $F(I)$ is Cohen-Macaulay if and only if $\m I^2 \cap JI = \m IJ$.
  \item If $\lambda(\m I/\m J) = 2$ and $\m
    I^2 = \m JI$, then $F(I)$ is Cohen-Macaulay.
\end{enumerate}
\end{proposition}

\begin{proof}
(1) By Proposition 5.4 of \cite{avj-jv1}, $F(I)$ is Cohen-Macaulay if and
only if $\lambda(\m I^n + JI^{n-1}/JI^{n-1}) = 1$ for all $n = 1,
\ldots, r^\m_J(I)$. Since $r^\m_J(I) = 2$, this equation translates to
$\lambda(\m I^2 + JI/JI) = 1$. Since $\lambda(\m I^2/\m JI) = 1$, this
is equivalent to $\m I^2 \cap JI = \m JI$.
\vskip 2mm
\noindent
(2) Note that in this case $r^\m_J(I) = 1$. By Corollary
\ref{cor-redn}(2), $F(I)$ is Cohen-Macaulay.
\end{proof}
\section{Examples}

\begin{example} (Example 4.6, \cite{avj-jv1})\label{ex1}
Let $R=\mathbb{Q}[\![x,y,z]\!]$. Let
$I=(-x^2+y^2,-y^2+z^2,xy,yz,zx)$ and $J=(-x^2+y^2,-y^2+z^2,xy)$. Since
$I^3=JI^2$, $J$ is a minimal reduction of $I$. Also $\lambda(\m I/\m
J)=1$ and $\m I^2=\m JI$. Therefore $\lambda(\m
I/\m J)=1$. It has been shown in \cite{avj-jv1} that $\depth F(I) =
1$ and $\depth G(I) = 0$. Here the $\m$-reduction number of $I$ is 1.
This example shows that even if the $\m$-reduction number of $I$ is 1,
the depth of fiber cone can be quite low without high depth assumption
on the associated graded ring.
\end{example}

\begin{example}\label{ex2}
Let $R = \mathbb{Q}[\![X,~Y,~Z,~W,~U]\!]/K$, where $K = (-X^3Z + Y^3,~ X^5
- Z^2,~ -X^2Y^3 + Z^3,~ -X^4Y^2 + ZW,~ -X^2Z^2 + YW,~ -Y^2Z + XW,~ XYZ^3 -
W^2,~ X^3YW - Z^4,~ Z^5 - Y^4W,~ -Y^5 + X^4W)$. Then it can be seen that
$R \cong \mathbb{Q}[\![t^6, ~t^{11},~ t^{15},~ t^{31},~ u]\!]$, is
Cohen-Macaulay of dimension $2$. Let $x = X + K,\; y = Y+K,~ z = Z +
K,~ w = W+K$ and $u = U + K$. Let $\m = (x,y,z,w,u), \; I = (x,y,w,u)$ and
$J = (x,u)$.  Then it can be seen that $I^3 = JI$, $\lambda(\m I/\m J)
= 2,\; \lambda(\m I^2 \cap J/\m JI) = 1, \; \m I^{n+1} \cap J = \m
JI^n$ for all $n \geq 2$ and $\lambda(\m I^2/\m JI) =
1$. It can also be seen that $I^2 \cap J = JI$. Therefore $I^n \cap J
= JI^{n-1}$ for all $n \geq 1$. This implies that $G(I)$ is
Cohen-Macaulay, i.e., $\depth G(I) = 2$. Since $\m I^2 \cap J \neq \m
JI$, $F(I)$ is not Cohen-Macaulay. It can also be verified that $\m
I^n \cap (u) \subseteq \m I^{n-1}J$ for all $n \geq 1$ so that $u^o$
is regular in $F(I)$. Hence $\depth F(I) = 1$.  Therefore $\sum_{n\geq
2} \lambda(\m I^n \cap J/\m JI^{n-1}) = 1$ with $\depth G(I) = 2$ and
$\depth F(I) = 1$. 
\end{example}

\begin{example}\label{ex4}
Let $R = \mathbb{Q}[\![x,y,z]\!], \; I = (x^3, y^3, z^3, xyz, xz^2, yz^2)$ and
$J = (x^3, y^3, z^3)$. Then it can be checked that $r_J(I) = r^\m_J(I)
= 2, \; \lambda(\m I^2 \cap J/\m JI) = 1$ and $\m I^n \cap J = \m
JI^{n-1}$ for all $n \geq 3$. It can also be verified that $I^2 \cap J
\neq JI$. Therefore $G(I)$ is not Cohen-Macaulay. It can also be seen
that the Hilbert series
\begin{eqnarray*}
  H(G(I),t) & = & (1-t)^{-1}H(G(I/(x^3)),t) \\
  & = & (1-t)^{-2}H(G(I/(x^3, y^3)),t) \\
  & = & \frac{15+6t+6t^2}{(1-t)^3}.
\end{eqnarray*}
Therefore $(x^3)^*, (y^3)^*$ is a regular sequence in $G(I)$. Hence
$\depth G(I) = 2$. Also, it can be computed that $e_0(F(I)) = 9 = 1 +
\lambda(I/J+\m I) + \lambda(I^2/JI + \m I^2)$. Therefore by Theorem
2.1 of \cite{drv}, $F(I)$ is Cohen-Macaulay.
\end{example}

\begin{example}\label{ex3}
Let $R = \mathbb{Q}[\![x,y]\!], \; I = (x^5, x^3y^2, x^2y^4, y^5)$ and $J =
(x^5, y^5)$. Then it can be verified that $I^5 = JI^4, \;
\lambda(\m I^2 \cap J/\m JI) = 1$ and $\m I^{n+1} = \m JI^n$ for all
$n \geq 2$. Let $u = x^5+y^5$. Then the Hilbert Series 
\begin{eqnarray*}
  H(G(I),t) & = & (1-t)^{-1}H(G(I/u),t) \\
& = & \frac{18 + 6t + t^4}{(1-t)^2}.
\end{eqnarray*}
Therefore $u^*$ is regular in $G(I)$. Since $I^2 \cap J \neq JI$, by
Valabrega-Valla condition, $G(I)$ is not Cohen-Macaulay. Hence $\depth
G(I) = 1$. Therefore by Theorem \ref{main}, $\depth F(I) \geq 1$. It
can also be verified that $\mu(I^3) = 11 < 13 = \sum_{n=0}^3
(n+1)\lambda(I^{3-n}/JI^{3-n-1}+\m I^{3-n})$. Therefore by Theorem 2.1
of \cite{drv}, $F(I)$ is not Cohen-Macaulay. Hence $\depth F(I) = 1$.
This shows that the lower bound for the depth of $F(I)$ given in
Theorem \ref{main} is sharp.
\end{example}

\end{document}